\documentclass[10pt]{article}
\usepackage{preambule}
\title{Remark on the semilinear ill-posedness for a periodic higher order KP-I equation}

\date{}
\author{Tristan Robert\footnote{Contact : tristan.robert@u-cergy.fr}\\ \emph{Université de Cergy-Pontoise}\\ \emph{Laboratoire AGM}\\ \emph{2 av. Adolphe Chauvin, 95302 Cergy-Pontoise Cedex, France}}

\begin{document}
\maketitle

\selectlanguage{english}
\begin{abstract}
We prove that, for some irrational torus, the flow map of the periodic fifth-order KP-I equation is not locally uniformly continuous on the energy space, even on the hyperplanes of fixed $x$-mean value.
\end{abstract}

\selectlanguage{french}
\begin{center}
\large{\emph{Remarque sur le caractère semi-linéairement mal posé pour une équation KP-I périodique d'ordre supérieur}}
\end{center}
\begin{abstract}
On montre que, pour un tore irrationnel bien choisi, le flot pour l'équation KP-I d'ordre 5 périodique n'est pas localement uniformément continu sur l'espace d'énergie, même sur les hyperplans de données initiales à moyenne en $x$ fixée.
\end{abstract}

\section{Introduction}
The study of well-posedness for nonlinear dispersive equations has seen constant progress during the last few decades. A cornerstone in the low regularity Cauchy theory for these equations has been the work of Bourgain \cite{bourgain1993kdv} who developped a remarkably effective method to prove local well-posedness, based on a fixed point argument in a function space tailored to the linear part of the equation. A consequence of this approach is to provide a flow map which is continuous on the Sobolev spaces, and even locally Lipschitz continuous. In that case we say that the problem is semilinearly well-posed (see \cite{Tzvetkov2004ill}).

In the early 2000's, though, some examples arose showing that this behaviour is not universal. First, the failure of $\mathcal{C}^2$ regularity for the KdV equation below $\dot{H}^{-3/4}$ \cite{Tzvetkov1999} (which was already known for $\mathcal{C}^3$ \cite{Bourgain1997}) suggested that not only the bilinear estimate in Bourgain spaces may not hold all the way down to the scaling critical regularity (which corresponds to a control on the second Picard iterate and thus on the level of regularity $\mathcal{C}^2$), but that the semilinear ill-posedness may appear. Then, the case of the periodic fifth-order KP-I equation \cite{SautTzvetkov2001} showed that the bilinear estimate in standard Bourgain spaces can even fail at any regularity. Of course, this is a priori not an objection to still recover a smooth flow map by performing an iteration procedure in other functional spaces : for the KP-II equation, the bilinear estimate is no longer true in standard Bourgain spaces modeled on the anisotropic Sobolev space $H^{s_1,s_2}(\R^2)$ for $s_1<-1/3$ \cite{Takaoka}, yet a Picard iteration can be performed in well chosen spaces for data in the scaling at the scaling critical rgularity $s_1=-1/2$ \cite{HadacHerrKoch}. Nevertheless, combining the two ideas above, Molinet, Saut and Tzvetkov  \cite{MST2002} proved that, as far as the KP-I equation is concerned, the failure of the bilinear estimate is not just of technical nature, since once again the flow map cannot be of class $\mathcal{C}^2$ in any Sobolev space. Then another counter-example was given by the Benjamin-Ono equation \cite{MolinetSautTzvetkov2001}. To treat both these equations, one is then forced to give up on the implementation of the contraction principle, which motivated the development of new methods based on compactness rather than completeness to attack those problems.

Consequently, the Cauchy problem for both these equations has been extensively investigated. The latter was shown to be globally well-posed in both $L^2(\R)$ \cite{IonescuKenig2007BO} and $L^2(\T)$ \cite{Molinet2008}. On the real line, its quasilinear behaviour (in the sense of lack of regularity of the flow map) was further examined in \cite{KochTzvetkov2003BO}, where it was shown that the transport effect due to the derivative in the nonlinearity leads to a change of speed of the plane waves, which in turn contradicts the local uniform continuity. Regarding the periodic Benjamin-ono equation, it was also proved in \cite{Molinet2008} that the flow map is not locally uniformly continuous on the  whole space $L^2(\T)$, yet it is on the subspace of zero mean value data. For the former, which is known to be globally well-posed in the energy spaces $\E(\R^2)$ \cite{IonescuKenigTataru2008} and $\E(\R\times\T)$ \cite{Article1} associated with its Hamiltonian structure, the same effect has been exploited in \cite{KochTzvetkov2008} (along with a  transverse effect) to show the lack of local uniform continuity of the flow map.

All these semilinear ill-posedness results thus rely on the failure of uniform continuity for the Galilean transformation
\begin{equation}\label{definition galilean transform}
\G_t^{\pm} : u_0 \in H^s(\T) \mapsto u_0\left(\cdot \pm t\fint_{\T}u_0(x)\dx\right) \mp \fint_{\T}u_0(x)\dx,
\end{equation}
which is well-defined for $s\supeg 0$.

Indeed, for $n\in\N^*$, take 
\[u_n(x) := n^{-s}\cos(nx) +n^{-1} \text{ and }v_n(x) := n^{-s}\cos(nx),\]
then $u_n,v_n$ are uniformly bounded in $H^s(\T)$, satisfy 
\[\normH{s}{u_n-v_n}\sim n^{-1}\underset{n\rightarrow+\infty}{\longrightarrow}0,\] but for $t\in ]0;1]$
\begin{multline*}
\normH{s}{\G_t^+(u_n)-\G_t^+(v_n)}=c\normH{s}{n^{-s}\left\{\cos\left(n(x+n^{-1}t)\right)-\cos(nx)\right\}}\\
\gtrsim |\sin(t)|>0
\end{multline*}
which shows the non-uniform continuity of $\G_t^+$ for $t\in ]0;1]$ (whereas for $t=0$ then $\G_0^{\pm}$ is Lipschitz continuous on $H^s(\T)$ for any $s\supeg 0$).

In particular, the general strategy to study the Cauchy problem is to construct a flow map $\Phi_t^0 : u_0 \mapsto u(t)$ on the subspace $L^2_0(\T)\subset L^2(\T)$ of zero mean value data, and then to obtain a flow map on the whole space via the formula
\begin{equation}\label{definition flot espace total}
\Phi_t := \G_t^-\circ \Phi_t^0\circ \G_0^+.
\end{equation}
This is possible since the mean value is a constant of the motion. Still, for the periodic Benajmin-Ono equation, the map $\Phi_t^0$ itself is actually Lipschitz continuous \cite{Molinet2008}. The argument above is in fact quite general : it applies for any periodic Hamiltonian equation under the form
\begin{equation}\label{equation hamiltonienne}
\drt u = \drx\nabla\H(u(t))
\end{equation}
given a Hamiltonian functional $\H : u_0(x,y)\in H^s(\T^{1+d}) \mapsto \H(u_0)\in\R$, for which the $x$-mean value is independent of $t$ (as can easily be seen by integrating (\ref{equation hamiltonienne}) in $x$). For example, this is the case for the periodic KdV equation and for the KP-II equation on the torus, which are nonetheless globally semilinearly well-posed on $L_0^2(\T)$ \cite{bourgain1993kdv} and $L_0^2(\T^2)$ \cite{bourgain1993kp} respectively. As for the BO and KP-I equations on $\R$ and $\R^2$ respectively, the failure of uniform continuity is obtained by using a properly localized version of the Galilean transformation above.

Note that this argument is very similar to the one observed by Herr \cite{Herr2006phd} concerning the derivative nonlinear Schrödinger equation, with the Galilean transformation $\G_t^{\pm}$ being replaced by the gauge transformation
\[\mathbf{G}_t^{\pm} := u_0\in H^s(\T) \mapsto \e^{\pm i\mathcal{I}(u)}u_0(\cdot \pm t\fint_{\T}|u|^2\dx),\]
with
\[\mathcal{I}(u) = \fint_{\T}\drx^{-1}\left(|u|^2-\fint_{\T}|u|^2\dx'\right)\dx\]
which leads to the failure of local uniform continuity of the flow map,  yet the local uniform continuity is recovered on the spheres of data with prescibed $L^2(\T)$ norm (which replace the above hyperplanes in this case).

Our purpose here is then to come back to the first known example of failure of the bilinear estimate, namely the periodic fifth-order KP-I equation
\begin{equation}\label{equation 5kp1}
\drt u -\drx^5u-\drx^{-1}\dry^2 u + u\drx u =0,~(t,x,y)\in\R\times\T^2.
\end{equation}
This model admits the Hamiltonian structure (\ref{equation hamiltonienne}) given by the Hamiltonian
\begin{equation}\label{definition hamiltonien}
\H(u_0) := \frac{1}{2}\normL{2}{\drx^2u_0}^2+\frac{1}{2}\normL{2}{\drx^{-1}\dry u_0}^2-\frac{1}{6}\int_{\T^2}u_0(x,y)^3\dx\dy,
\end{equation}
where the operator $\drx^{-1}\dry$ is well-defined on
\[\D_0'(\T^2):=\left\{u_0\in\mathcal{D}'(\T^2),~\widehat{u_0}(0,n)=0~\forall n\neq 0\right\}\]
as the Fourier multiplier with symbol $n/m$.

The local well-posedness for this equation was first studied in \cite{IorioNunes} for data in $H^s(\T^2)\cap \mathcal{D}_0'(\T^2)$ for $s>2$. In \cite{Article2}, we constructed a global flow map on the energy space 
\begin{equation}\label{definition espace energie}
\E^2(\T^2) := \left\{u_0\in L^2(\T^2)\cap \D_0'(\T^2),~\H(u_0)<+\infty\right\},
\end{equation}
endowed with the norm 
\begin{equation}\label{definition espace energy}
\norme{\E^2}{u_0}^2:=\normL{2}{u}^2+\normL{2}{\drx^2 u_0}^2 + \normL{2}{\drx^{-1}\dry u_0}^2,
\end{equation}
and proved persistence of regularity in the Banach scale $\E^{\sigma}$, $\sigma\supeg 2$ of functions with finite norm
\begin{equation}\label{definition E sigma}
\norme{\E^{\sigma}}{u_0}^2:=\normL{2}{u}^2+\normL{2}{\drx^{\sigma} u_0}^2 + \normL{2}{\drx^{-1}\dry u_0}^2+\normL{2}{\drx^{\sigma-3}\dry u_0}^2.
\end{equation}
This flow map was constructed on the whole energy space by the procedure described above (\ref{definition flot espace total}) and as such is indeed not uniformly continuous on $\E^{\sigma}(\T^2)$. Note that, from the definition of $\mathcal{D}_0'(\T^2)$ and the Hamiltonian structure (\ref{equation hamiltonienne}), the $x$-mean value is actually a constant of the motion, independent of both $t$ and $y$, thus $\G_t^{\pm}$ is well-defined. The aim of this note is to show that the quasilinear behaviour of equation (\ref{equation 5kp1}) is actually more involved, by proving the failure of uniform continuity of the flow map $\Phi_t^0$ defined on the hyperplane $\E^{\sigma}_0(\T^2)\subset \E^{\sigma}(\T^2)$, $\sigma\supeg 2$, of zero $x$-mean value data :
\begin{theoreme}\label{theorem}
There exists $\lambda>0$ such that for any $\sigma\supeg 2$, there exists two positive constants $c$ and $C$ and two sequences $(u_n)$ and $(v_n)$ of solutions to (\ref{equation 5kp1}) in $\mathcal{C}([0;1],\E^{\sigma}_0(\Tl^2))$ such that
\begin{equation}\label{estimation borne uniforme solutions}
\sup_{t\in[0;1]}\norme{\E^{\sigma}}{u_n(t)}+\norme{\E^{\sigma}}{v_n(t)}\infeg C,
\end{equation}
and satisfying initially
\begin{equation}\label{estimation solutions initiales}
\lim_{n\rightarrow +\infty}\norme{\E^{\sigma}}{u_n(0)-v_n(0)}=0,
\end{equation}
but such that for every $t\in[0;1]$,
\begin{equation}\label{estimation failure uniform continuity}
\liminf_{n\rightarrow +\infty}\norme{\E^{\sigma}}{u_n(t)-v_n(t)}\supeg c|t|.
\end{equation}
\end{theoreme}
Here, we will work on an irrational torus $\Tl^2 := \T\times\lambda^{-1}\T$ for some $\lambda>0$. Note that the construction of the flow map in \cite{Article2} is performed for a square torus but is completely insensitive to the choice of the periods of the initial data, thus it can be adapted on $\Tl^2$ in a straightforward manner.

\section{Outline of the proof}
Let us now discuss the strategy of the proof of theorem~\ref{theorem}. As explained above, working on $\E_0^{\sigma}(\Tl^2)$ rules out the transport effect due to the Burgers type nonlinearity. Here, the main nonlinear phenomenon comes from the resonant low (but non zero)-high frequency interaction.

Indeed, for the KP-II equation, the proof of global semilinear well-posedness of \cite{bourgain1993kp} heavily uses an algebraic feature of the equation : the symbol $\widetilde{\omega}(m,n)=m^3-n^2/m$ satisfies the nonresonant relation
\begin{multline}\label{definition Omega KPII}
\left|\widetilde{\Omega}(m_1,n_1,m_2,n_2)\right| := \left|\widetilde{\omega}(m_1+m_2,n_1+n_2)-\widetilde{\omega}(m_1,n_1)-\widetilde{\omega}(m_2,n_2)\right|\\
= \left|\frac{m_1m_2}{m_1+m_2}\left\{3(m_1+m_2)^2+\left(\frac{n_1}{m_1}-\frac{n_2}{m_2}\right)^2\right\}\right|\\
\gtrsim \left|m_1m_2(m_1+m_2)\right|,
\end{multline}
which provides a smoothing effect in the nonlinear interaction which compensates for the derivative loss.

For equation (\ref{equation 5kp1}), the symbol reads 
\begin{equation}\label{definition symbole 5kp1}
\omega(m,n)=m^5+\frac{n^2}{m},
\end{equation}
so that the resonant function becomes
\begin{multline}\label{definition fonction resonance}
\Omega(m_1,n_1,m_2,n_2) \\= \frac{m_1m_2}{m_1+m_2}\left\{5(m_1+m_2)^2(m_1^2+m_1m_2+m_2^2) - \left(\frac{n_1}{m_1}-\frac{n_2}{m_2}\right)^2\right\},
\end{multline}
which now enjoys a large set of resonant frequencies $(m_1,n_1,m_2,n_2)$ which annul $\Omega$. In particular, for $n\in\N$ we can choose $\alpha(n)\in\Zl$ such that
\begin{equation}\label{definition alpha}
\Omega(1,0,n,\alpha(n))=0.
\end{equation}
This means that the nonlinear interaction between the linear solution with frequency $(1,0)$ and the one with $(n,\alpha(n))$ produces a linear solution with frequency $(n+1,\alpha(n))$. This particular interaction was already exploited in \cite{SautTzvetkov2001} to prove the failure of the bilinear estimates in Bourgain spaces. To get the failure of local uniform continuity, we will analyze more closely how the initial data considered in \cite{SautTzvetkov2001} evolves under the nonlinear flow constructed in \cite{Article2}. 

More precisely, in section~\ref{section construction} we construct a family of functions which agree at time zero with the initial data given by the two modes considered above, and who solve the equation (\ref{equation 5kp1}) up to a sufficiently small error. The ansatz for this construction is to compute the first Picard iterates. Of course, the argument in \cite{SautTzvetkov2001} shows that this iteration scheme does not converge, yet the analysis in \cite{Article2} relies on this iteration on small times of order $O(n^{-2})$. Here, in order to have a good approximation up to time $O(1)$ we slightly damp the low frequency component. The first iterate is simply the linear evolution, and the second iterate describes the nonlinear interaction between the linear solutions of each frequency, in particular the low-high frequency interaction produces a new linear solution amplified linearly in time (which is the reason for the divergence of the iteration scheme). We will actually take for the low frequency component the genuine solution emanating from the low frequency mode to kill the low-low interaction, which remains of the same order as the main nonlinear effect. For complex valued solutions, we can stop the approximation here, but in order to work with real valued solutions, we modulate the high frequency mode with an oscillating term in time to absorb some error terms that appear with this latter constraint, and we will also take into account the main component of the third iterate, in order to have an appropriate error after plugging this approximate solution in the equation. In section~\ref{section comparaison}, we then compare these functions with the nonlinear flow of \cite{Article2} applied to the same data by using a standard energy method, and then we conclude the proof of theorem~\ref{theorem} in section~\ref{section preuve}.
\section{Construction of a family of approximate solutions}\label{section construction}
Let us fix $\sigma\supeg 2$. For $\theta \in [-1;1]$ and $n\in\N^*$, let us define the family of functions on $[0;1]\times\Tl^2$ by
\begin{multline}\label{definition solution approchee}
u_{\theta,n}(t,x,y):= \Phi_t\left[\theta n^{-1}\cos(x)\right] + \cos\left(\frac{\theta}{2}t\right)n^{-\sigma}\cos\left(\varphi_{n}(t,x,y)\right)\\+\sin\left(\frac{\theta}{2}t\right)n^{-\sigma}\sin\left(\varphi_{n+1}(t,x,y)\right)+R_{\theta,n}(t,x,y),
\end{multline}
where the phase functions are given by
\[\varphi_1 := x+t,~\varphi_{n} := nx+\alpha(n)y+\omega(n,\alpha(n))t\]
and
\[\varphi_{n+1} :=(n+1)x + \alpha(n)y+\omega(n+1,\alpha(n))t.\]
They are the the phase functions of linear solutions, that is $\cos(\varphi_k)$ solves the equation
\[(\drt -\L)\cos\varphi_k = 0,\]
where $\L = -\drx^5 - \drx^{-1}\dry^2$ is the linear operator in (\ref{equation 5kp1}). 
%Indeed, for $v\in \E^{\infty}(\Tl^2)\cap\D_0'(\Tl^2)$ (see the definition below), we have the equivalent representation for the operator $\drx^{-1}$ as the unique primitive which has zero $x$ mean value.

The rest is given by
\begin{equation}\label{definition reste}
R_{\theta,n} =n^{-\sigma}\left\{\cos\left(\frac{\theta}{2}t\right)\Omega_{n-1}^{-1}\cos\left(\varphi_n-\varphi_1\right) + \sin\left(\frac{\theta}{2}t\right)\Omega_{n+1}^{-1}\sin\left(\varphi_{n+1}+\varphi_1\right)\right\},
\end{equation}
where
\[\Omega_{n\pm 1}  := \pm\Omega(1,0,n\pm 1,\alpha(n)).\]
Note that a straightforward computation (see below for the definition of $\alpha(n)$) gives
\begin{equation}\label{estimation Omega npm1}
\left|\Omega_{n\pm 1}\right|\sim n^3,
\end{equation}
provided $\alpha(n)$ satisfies (\ref{definition alpha}). Of course in this case we also have
\begin{equation}\label{equation phases}
\varphi_1+\varphi_{n}=\varphi_{n+1}.
\end{equation}
\subsection{On the choice of the period}
 In order to annul the resonnant function, we have the ansatz \[\alpha(n)=n(n+1)\sqrt{5n^2+5n+5}\in\Zl\]
We are thus looking for a $\lambda>0$ such that for $n\in\N$ then $\sqrt{5n^2+5n+5} = \lambda n_1$ with $n_1\in\Z$. If we take $\lambda = \sqrt{5\ell}$ with $\ell\in\N$, setting then $X = 2n+1$ and $Y = 2n_1$ we are thus left with finding the integer solutions to
\begin{equation}\label{equation hyperbole}
X^2-\ell Y^2 = -3
\end{equation}
 Note that we want $n\rightarrow +\infty$ in the following, so we have to choose $\ell\in\N$ such that the above hyperbola  has an infinite number of integer points. Now, it is well known that the solutions $(X_k,Y_k)_{k\in\N}$ to (\ref{equation hyperbole}) are given by 
\[X_k+Y_k\sqrt{\ell} = (Y_0+\sqrt{\ell}X_0)(u_k+\sqrt{\ell}v_k),\]
where $(X_0,Y_0)$ is a particular solution and $(u_k,v_k)$ is a solution to Pell's equation $u^2-\ell v^2=1$. If $\ell$ is square free, Pell's equation has an infinite number of solutions $u_k+\sqrt{\ell}v_k = (u_0+\sqrt{\ell}v_0)^k$ for all $k\in\N$, where $(u_0,v_0)$ is the fundamental solution, thus it is enough to find a square free integer $\ell$ such that (\ref{equation hyperbole}) has at least one solution. For example, we can take $\ell = 7$ and $(X_0,Y_0)=(2,1)$. This choice of $\lambda = \sqrt{35}$ provides an infinite set of numbers $\{N_k\}\subset \N^{\N}$ such that
\begin{equation}
\alpha(N_k):= N_k(N_k+1)\sqrt{5N_k^2+5N_k+5}\in \lambda\Z.
\end{equation}
Thus in the following we will work with the functions $u_{\theta,n}$ (\ref{definition solution approchee}) for $n\in \{N_k\}$.
\subsection{Estimates on the approximate solutions}
First, let us recall the precise statement of the definition of the nonlinear flow \cite[Theorem 1.1 (a), Proposition 6.2]{Article2} :
\newpage
\begin{theoreme}\label{theoreme flot}For any ${\displaystyle u_0\in\E^{\infty}(\Tl^2):= \bigcap_{\sigma \supeg 2}\E^{\sigma}(\Tl^2)}$, there exists a unique global smooth solution \[u=:\Phi^{\infty}(u_0)\in\mathcal{C}(\R,\E^{\infty}(\Tl^2))\] to (\ref{equation 5kp1}) and moreover there exists a positive 
\begin{equation}\label{estimation temps existence}
T=T(\norme{\E^2}{u_0})\sim \crochet{\norme{\E^2}{u_0}}^{-\mu}
\end{equation}
for some $\mu>0$ such that for any $\sigma\supeg 2$ we have
\begin{equation}\label{estimation energie avec donnee}
\norme{L^{\infty}_T\E^{\sigma}}{\Phi^{\infty}(u_0)}\infeg C_{\sigma}\norme{\E^{\sigma}}{u_0}.
\end{equation}
\end{theoreme}
Next, we prove several bounds on $u_{\theta,n}$. First, as explained above, for the low frequency part we took the nonlinear solution instead of the linear one to avoid the contribution of the low-low interaction. The next lemma shows that in all the other nonlinear interactions we can replace the former by the latter up to a manageable error :
\begin{lemme}\label{lemme estimation basse frequence}
The nonlinear low frequency solution $u_1 := \Phi_t\left[\theta n^{-1}\cos(x)\right]$ is close enough to the linear solution $\widetilde{u_1} := \theta n^{-1} \cos\varphi_1$, namely
\begin{equation}\label{estimation erreur basse frequence nonlineaire}
\normL{2}{u_1-\widetilde{u_1}} \lesssim n^{-2},
\end{equation}
uniformly in $\theta\in[-1;1]$ and $t\in [0;1]$.
\end{lemme}
\begin{proof}
Let us write $v_1 := u_1-\widetilde{u_1}$, then $v_1$ solves the equation
\[(\drt - \L)v_1 = -u_1\drx u_1.\]
Thus, multiplying by $v_1$, integrating, using the skew-symmetry of $\L$ and Cauchy-Schwarz inequality, we get
\[\ddt \normL{2}{v_1(t)}^2\lesssim \normL{2}{v_1(t)}\normL{2}{u_1(t)}\normL{\infty}{\drx u_1(t)}.\]
The $L^2$ norm of $u_1$ is preserved by the flow, thus this term is $O(n^{-1})$. Next, using a Sobolev inequality and (\ref{estimation energie avec donnee}), the last term above can be estimated as
\[\norme{L^1([0;1])L^{\infty}}{\drx u_1(t)}\lesssim \norme{L^{\infty}([0;1])\E^{10}}{u_1}\lesssim \norme{\E^{10}}{u_1(0)}=O(n^{-1})\]
Using these bounds and integrating the first estimate on $[0;1]$ and using that $v_1(0)=0$, we finally get the bound
\[\norme{L^{\infty}([0;1]L^2)}{v_1}^2 \lesssim n^{-2}\norme{L^1([0;1]L^2)}{v_1}\] which provides (\ref{estimation erreur basse frequence nonlineaire}).
\end{proof}
Once we control the low frequency part, we can obtain the main bounds for the approximate solutions $u_{\theta,n}$ :
\begin{lemme}\label{lemme estimation uapp}
For $u_{\theta,n}$ defined in (\ref{definition solution approchee}), the following estimates hold uniformly in $\theta\in [-1;1]$, $t\in[0;1]$ and $n\in \{N_k\}$ :
\begin{equation}\label{estimation norme solution approchee}
\norme{\E^{\sigma}}{u_{\theta,n}}\lesssim 1
\end{equation}
and
\begin{equation}\label{estimation erreur equation}
\normL{2}{(\drt - \L)u_{\theta,n}+u_{\theta,n}\drx u_{\theta,n}}\lesssim n^{-\sigma-1}.
\end{equation}
\end{lemme}
\begin{proof}
The only nontrivial fact in (\ref{estimation norme solution approchee}) is the bound on the low frequency part $u_{\theta,n}$, which is again a consequence of (\ref{estimation energie avec donnee}). So it remains to prove the main estimate (\ref{estimation erreur equation}).

Let $u_i$, $i=1,2,3$ be the modes in (\ref{definition solution approchee}) and $\widetilde{u_1}$ be as in the previous lemma. Then the term in the left-hand side of (\ref{estimation erreur equation}) can be written as
\begin{multline}
(\drt - \L)u_{\theta,n} + u_{\theta,n}\drx u_{\theta,n}\\
= (\drt - \L)(u_2+u_3+R_{\theta,n}) + \widetilde{u_1}\drx (u_2+u_3) + F_1+F_2,
\end{multline}
where we have used that $u_1$ is a genuine nonlinear solution, so that here the remainder terms are
\[F_1 = (u_1-\widetilde{u_1})\drx (u_2+u_3 + R_{\theta,n})
\text{ and }F_2 = (u_2+u_3+R_{\theta,n})\drx u_{\theta,n}.\]
Here and in the sequel, we write $O(n^{-\beta})$ various terms having an $L^2$ norm bounded by a constant (uniform in $\theta,t$ and $n$) times $n^{-\beta}$. 

The definition of $u_2$, $u_3$ and $R_{\theta,n}$ along with (\ref{estimation energie avec donnee}) for the low frequency term provide $F_2 = O(n^{-\sigma-1})$, and from lemma~\ref{lemme estimation basse frequence} we also have $F_1=O(n^{-\sigma-1})$, so it remains to compute the main contribution.

Since $\cos\varphi_k$ is a linear solution, we have
\begin{equation}\label{calcul u2 lineaire}
(\drt -\L) u_2 = -\frac{\theta}{2}\sin\left(\frac{\theta}{2}t\right)n^{-\sigma}\cos\varphi_n
\end{equation}
and
\begin{equation}\label{calcul u3 lineaire}
(\drt - \L)u_3 = \frac{\theta}{2}\cos\left(\frac{\theta}{2}t\right)n^{-\sigma}\sin\varphi_{n+1}.
\end{equation}
For the linear evolution of $R_{\theta,n}$, first note that the argument in the cosine and sine functions are respectively
\[\varphi_n-\varphi_1 = (n-1)x+\alpha(n)y +[\omega(n,\alpha(n))-\omega(1,0)]t\]
and
\[\varphi_{n+1}+\varphi_1 = (n+2)x+\alpha(n)y +[\omega(n+1,\alpha(n))+\omega(1,0)]t,\]
so that 
\begin{multline*}
(\drt -\L)\cos\left(\varphi_n-\varphi_1\right)=-[\omega(n,\alpha(n))-\omega(1,0) - \omega(n-1,\alpha(n))]\sin\left(\varphi_n-\varphi_1\right)\\
= \Omega_{n-1}\sin(\varphi_n-\varphi_1),
\end{multline*}
and similarly for the sine term. Thus
\begin{multline}\label{calcul R lineaire}
(\drt -\L)R_{\theta,n} = n^{-\sigma} \left\{\cos\left(\frac{\theta}{2}t\right)\sin\left(\varphi_n-\varphi_1\right)+\sin\left(\frac{\theta}{2}t\right)\cos\left(\varphi_{n+1}+\varphi_1\right)\right\}\\+O(n^{-\sigma-3})
\end{multline}
thanks to (\ref{estimation Omega npm1}).

Next, the nonlinear interactions are
\begin{multline}\label{calcul u2 nonlineaire}
\widetilde{u_1}\drx u_2 = -\theta \cos\left(\frac{\theta}{2}t\right) n^{-\sigma}\cos\varphi_1\cdot\sin \varphi_n\\
= -\frac{\theta}{2}\cos\left(\frac{\theta}{2}t\right) n^{-\sigma}\sin \varphi_{n+1}-\frac{\theta}{2}\cos\left(\frac{\theta}{2}t\right) n^{-\sigma}\sin\left(\varphi_{n}-\varphi_1\right),
\end{multline}
where we have used (\ref{equation phases}). Similarly,
\begin{equation}\label{calcul u3 nonlineaire}
\widetilde{u_1}\drx u_3 = \frac{\theta}{2}\sin\left(\frac{\theta}{2}t\right) n^{-\sigma}\cos\left(\varphi_{n+1}+\varphi_1\right) +\frac{\theta}{2}\sin\left(\frac{\theta}{2}t\right) n^{-\sigma}\cos\left(\varphi_{n}\right) +O(n^{-\sigma-1}).
\end{equation}
Thus we see that (\ref{calcul u2 lineaire})-(\ref{calcul u3 lineaire}) annul the main nonlinear terms in (\ref{calcul u2 nonlineaire})-(\ref{calcul u3 nonlineaire}), whereas $R_{\theta,n}$ (\ref{calcul R lineaire}) deals with the remainders in these interactions.

Finaly, summing (\ref{calcul u2 lineaire})-(\ref{calcul u3 nonlineaire}) yields (\ref{estimation erreur equation}).
\end{proof}

\section{Comparison with the genuine solutions}\label{section comparaison}
Let $u := \Phi_t\left[u_{\theta,n}(0)\right]$ be the genuine solution arising from the same initial data as for $u_{\theta,n}$. Note that, since ${\displaystyle \norme{\E^2}{u_{\theta,n}(0)} \lesssim n^{2-\sigma}}$, $u$ is defined on the time interval $[0;1]$ thanks to (\ref{estimation temps existence}).
\begin{lemme}\label{lemme comparaison}
Let $v:= u- u_{\theta,n}$ be the difference between the genuine and the approximate solutions. Then there exists $\delta>0$ such that
\begin{equation}\label{estimation difference solutions}
\norme{L^{\infty}([0;1])L^2}{\drx^{\sigma} v}\lesssim n^{-\delta}
\end{equation}
\end{lemme}
\begin{proof}
First, from the definition of $u_{\theta,n}$ and (\ref{estimation energie avec donnee}) we have
\begin{multline}\label{estimation v}
\norme{L^{\infty}([0;1])L^2}{\drx^{\sigma+1} v}\infeg \norme{L^{\infty}([0;1])\E^{\sigma+1}}{u(t)}+O(n)\\
\lesssim \norme{\E^{\sigma+1}}{u(0)}+O(n)=O(n),
\end{multline}
uniformly in $\theta\in[-1;1]$. Moreover, $v$ solves the equation
\[(\drt - \L)v + v\drx v +\drx\left(u_{\theta,n}\cdot v\right) + G =0,\]
where
\[G := (\drt -\L)u_{\theta,n}+u_{\theta,n}\drx u_{\theta,n}.\] 
Thus, multiplying by $v$ and integrating over $\Tl^2$, we get from the skew-symmetry of $\L$ the energy bound
\[\ddt \normL{2}{v(t)}^2 \lesssim \langle v,v\drx v\rangle + \langle v,\drx (u_{\theta,n}\cdot v)\rangle+\langle v,G\rangle,\]
where $\langle \cdot ,\cdot\rangle$ stands for the scalar product in $L^2(\Tl^2)$.

By integrating by parts, the first term is zero and the second one is bounded by
\[\normL{\infty}{\drx u_{\theta,n}(t)}\normL{2}{v(t)}^2 \lesssim  n^{1-\sigma}\normL{2}{v(t)}^2.\]
Thus, using Gronwall's inequality, we get
\[\norme{L^{\infty}([0;1],L^2)}{v}\lesssim \norme{L^{1}([0;1],L^2)}{G}\lesssim n^{-\sigma-1}\]
thanks to lemma~\ref{lemme estimation uapp}. Finally, (\ref{estimation difference solutions}) follows from interpolating the last bound with (\ref{estimation v}).
\end{proof}

\section{Proof of the main theorem}\label{section preuve}
Let us now conclude the proof of theorem~\ref{theorem}. Define $u_n := \Phi_t[u_{-1,n}(0)]$ and $v_n:=\Phi_t[u_{1,n}(0)]$. Then (\ref{estimation borne uniforme solutions}) follows from (\ref{estimation energie avec donnee}), and (\ref{estimation solutions initiales}) from
\[\norme{\E^{\sigma}}{u_n(0)-v_n(0)}= 2n^{-1}\norme{\E^{\sigma}}{\cos(x)}\rightarrow 0,\]
so it remains to prove (\ref{estimation failure uniform continuity}). We have
\[\norme{\E^{\sigma}}{u_n(t)-v_n(t)}\supeg \normL{2}{\drx^{\sigma}\left(u_n(t)-v_n(t)\right)}
\supeg \normL{2}{\drx^{\sigma} \left(u_{-1,n}(t)-u_{1,n}(t)\right)} -cn^{-\delta}\]
in view of the previous lemma. Using lemma~\ref{lemme estimation basse frequence} and the definition of $u_{\theta,n}$, this last term can be estimated from below by
\[O(n^{-1})+2\left|\sin(t/2)\right|\normL{2}{\sin^{(\sigma)}\varphi_{n+1}}+O(n^{-\delta})\gtrsim |t|+o(1),\]
which concludes the proof of (\ref{estimation failure uniform continuity}).

\section*{Aknowledgements}
The author is thankful to Nikolay Tzvetkov for his valuable comments.

\bibliographystyle{siam}
\bibliography{biblio}

\end{document}